\documentclass[a4paper]{article}
\usepackage[english]{babel}
\usepackage[dvips]{graphics,epsfig}
\usepackage{amsmath,graphicx}
\usepackage{amssymb,latexsym}
\usepackage{amsfonts}

\begin{document}
\begin{center}
\noindent {\Large\bf The Seven Classes of the Einstein Equations \\[1mm]
}

\end{center}
\bigskip

\noindent{\small{\bf Sergey E. Stepanov}}

\noindent{\small\it Professor, Chair of Mathematics, Finance
Academy under the Government of the Russian Federation, 49,
Leningradsky Prospect, Moscow, 125993 Russia.\\
 e-mail: stepanov@vtsnet.ru}
\medskip

\noindent{\small\bf Irina I. Tsyganok}

\noindent{\small\it Professor, Chair of General Scientific
Disciplines, Vladimir branch of Russian University of Cooperation,
14, Vorovskogo str., Vladimir, 6000001 Russia.\\
e-mail: i.i.tsyganok@mail.ru}

\medskip

\noindent{\small{\bf Abstract}. In current paper we refer to the
geometrical classification of the Einstein equations which has
been developed by one of the authors of this paper. This
classification was based on the classical theory for decomposition
of the tensor product of representations into irreducible
components, which is studied in the elementary representation
theory for orthogonal groups.
 We return to this result for more detailed investigation of classes of the Einstein equations.
\medskip

\medskip

}


\bigskip

\section*{ Introduction}

Method of classification of Riemannian manifolds, endowed with
tensor structure, on the basis of decomposition of irreducible
tensor product of representations (irreducible with respect to the
action of orthogonal group) into irreducible components has become
traditional in differential geometry (see, e.g., [1] - [5]).

This method is used in theoretical physics (see, for example, [6];
[7]; [8]). For example, in the paper [7] Einstein-Cartan manifolds
were classified on the basis of irreducible decomposition of the
torsion tensor of an affine-metric connection (irreducible with
respect to the action of the Lorentz group). Current results in
this modern today research direction in the General Relativity
Theory were systemized in paper [7].

At current paper we return back to the work [8] of one of the
authors and study in details proposed in this work classification
of The Einstein equations on the basis of method mentioned above.

\section {The Einstein equations}

 The space-time $(M,g)$  in the
General Relativity Theory is a smooth four-dimensional manifold
$M$, endowed with metrics $g$ of the Lorentz signature (- + + +).

In an arbitrary map  $(U,\varphi)$ with the local coordinate
system $\{x^0,x^1,x^2,x^3\}$ the metric $g$ of $(M,g)$  is defined
by its components $g_{ij}=g(\partial_i,\partial_j)$, where
$\partial_i=\partial/\partial x^i$ for $i,j,k,l, \ldots =0,1,2,3$.
These components define the Christoffel symbols
$\Gamma_{ij}^k=2^{-1}g^{kl}(\partial_i
g_{lj}+\partial_jg_{li}-\partial_l g_{ij})$ of the Levi-Civita
connection $\nabla$ and the operator of the covariant differention
defined by  $\nabla_i X^k=\partial_i X^k+\Gamma_{il}^jX^l$ for
$X=X^k\partial_k$.

In the General Relativity Theory a metric tensor $g$ is
interpreted as gravitational potential and is related by the {\it
Einstein equations}
$$
Ric-\frac{1}{2}sg=T:\quad R_{ij}-\frac{1}{2}s\;g_{ij}=T_{ij}\eqno
                                       (1.1)
$$
to mass-energy distribution, which generates the gravitational
field. Here $Ric$  is the Ricci tensor of metric $g$ of a
space-time $(M,g)$, which can be found from identity
$Ric(\partial_i,\partial_j)=R_{ij}=R^k_{ikj}$   for components
$R^l_{ikj}$ of torsion tensor $R$, which, in their turn, can be
found from equality $R^k_{lij}X^l=\nabla_i\nabla_j
X^k-\nabla_j\nabla_iX^k$, $s:=g^{ij}R_{ij}$ -- is scalar curvature
of metrics $g$ for $(g^{ij})=(g_{ij})^{-1}$ and, finally,
$T_{ij}=T(\partial_i,\partial_j)$ -- are components of known {\it
energy-momentum tensor $T$ of matter}.

 The curvature tensor $R$
satisfies well-known the Bianchi identities
 $$
 dR=0:\quad \nabla_mR^k_{lij}+\nabla_iR^k_{ljm}+\nabla_jR^k_{lmi}=0. \eqno   (1.2)
 $$
From these identities, in particular, following equalities can be
found
$$
d^*R=-d\;Ric: \quad
-\nabla_kR^k_{lij}=\nabla_jR_{li}-\nabla_iR_{lj}=0,\eqno (1.3)
$$
which lead to other the Bianchi identities
$$
2d^*Ric=-ds:\quad 2g^{kj}\nabla_jR_{ki}=\nabla_is. \eqno(1.4)
$$
Equations (1.1) are complemented by "conservation laws"
 $$
 d^*T=0:\quad -g^{ik}\nabla_iT_{kj}=0, \eqno(1.5)
 $$
which are derived from The Einstein equations on the basis of the
Bianchi identities (1.4).

\section{Invariantly defined seven classes \\of the
Einstein equations}

 In the paper [8] was proved that on a pseudo-Riemannian manifold $(M,g)$
   the bundle  $\Omega(M)\subset T^*M\otimes S^2M$, which fiber in each
    point $x\in M$  consists of linear mappings $\Omega:T_xM\rightarrow \mathbf{R}$
    such that $\Omega(X,Y,Z)=\Omega(X,Z,Y)$  and $\sum\limits_{i=1}^n
    \Omega(e_i,e_i,X)=0$  for arbitrary vectors $X,Y$  and $Z$ orthonormal
     basis $\{e_1,e_2,\dots,e_n\}$  of space  $T_xM$, has point-wise
     irreducible decomposition (irreducible with respect to the action of pseudo-orthogonal group)
     $\Omega(M)=\Omega_1(M)\oplus\Omega_2(M)\oplus\Omega_3(M)$.

If as $(M,g)$  we take space-time, then $\nabla T \in \Omega(M)$.
As the result six classes of the Einstein equations can be
determined via invariant approach $\Omega_\alpha$ and
$\Omega_\beta\oplus\Omega_\gamma$ ($\alpha,\beta,\gamma=1,2,3\mbox
{ и } \beta<\gamma$), for each of them covariant derivatives
$\nabla T$ of energy-momentum tensors $T$ of matter are
cross-section of relevant invariant subbundles
$\Omega_1(M),\Omega_2(M),\Omega_3(M)$ and their direct sums
$\Omega_1(M)\oplus\Omega_2(M)$, $\Omega_1(M)\oplus\Omega_3(M)$ и
$\Omega_2(M)\oplus\Omega_3(M)$. Seventh classes determined in [8]
satisfy the condition $\nabla T=0$.

\section{The class $\Omega_1$  and integrals of\\
geodesic equations }

Class $\Omega_1$ of the Einstein equations is selected via
condition (see, [8])
$$
\delta^*T=0:\quad \nabla_kT_{ij}+\nabla_iT_{jk}+\nabla_jT_{ki}=0.
\eqno(3.1)
$$
As  $trace_gT=-s$, then from (3.1) follows that  $d(trace_g
T)=2d^*T=0$, i.e. $s=const$. In this case equations (3.1) can be
re-written in the following form $\delta^*Ric=0: \nabla_kR_{ij}+
\nabla_iR_{jk}+\nabla_jR_{ki}=0$. Reverse is evident, i.e. from
condition $\delta^*Ric=0$ can be deduced that $s=const$, and then
equations (3.1).

If each solution $x^k=x^k(s)$ of geodesic equations on
pseudo-Riemannian manifold $(M,g)$ satisfies condition
$a_{ij}\frac{dx^i}{ds}\frac{dx^j}{ds}=const$ for symmetric tensor
field $a$ with local components $a_{ij}$, then it's said that such
geodesic equations admit {\it first quadratic integral} (see, for
example, [9]. In fact, following identity needs to be true
$\delta^*a=0$. At the same time tensor field  $a$, determined by
identity above, is called the {\it Killing tensor} (see [10]). In
our case, first quadratic integral of geodesic equations will be
$R_{ij}\frac{dx^i}{ds}\frac{dx^j}{ds}=const$ and Ricci tensor will
be consequently the Killing tensor.

 Following theorem was proved:

{\bf Theorem 1.} {\it The Einstein equations belong to the class
$\Omega_2$ if and only if in the space-time $(M,g)$ geodesic lines
equations admit following first quadratic integral
$R_{ij}\frac{dx^i}{ds}\frac{dx^j}{ds}=const$ for Ricci tensor
$R_{ij}=Ric(\partial_i,\partial_j)$ of metrics $g$.}

 {\bf Remark.} Let us
underline that theory of first integrals equations of geodesic and
symmetric tensor Killing fields has many applications in
mechanics, general relativity theory and other physics research
(see, for example, [10] and [11]).

\section{Class  $\Omega_2$ and the Young-Mills equations}

 Class $\Omega_2$  of The Einstein equations is selected via
  condition (see, [8])
$$
dT=0 :\quad \nabla_kT_{ij}=\nabla_iT_{kj}. \eqno(4.1)
 $$
Considering equations (3.1) together with The Einstein equations
we get that  $s=const$, and consequently we get following
equations
$$
d \;Ric=0:\quad \nabla_kR_{ij}=\nabla_iR_{kj}. \eqno(4.2)
$$
Equations (4.2) are known as {\it Codazzi equations} (see, for
details, [12]). Easy to prove that equations (4.1) follow from
Codazzi equations (4.2). In fact, from equations (4.2) follows
that  $d^*Ric=-ds$. Comparing this equation with Bianchi
identities (1.4) we conclude that $s=const$.  As the result
equations (4.1) will be true. So, equations (4.1) and (4.2) are
equivalent.

For constructing an example of space-time $(M,g)$  with $\nabla T
\in \Omega_2(M)$, we will take $n$-dimensional $(n\geq 4)$
conformally flat pseudo-Riemannian manifold  $(M,g)$, for which as
known (see, for example, [9]) following equations are fulfilled
$$
d\left[Ric-\frac{1}{2(n-1)}s\cdot g\right]=0.
$$
If we suppose here that $s=const$  we will get equations (4.2).
Therefore, on conformally flat space-time  $(M,g)$ with the
constant scalar curvature $s$ the Einstein equations belong to the
class  $\Omega_2$.

On arbitrary $n$-dimensional pseudo-Riemannian manifold $(M,g)$
the {\it Weyl projective curvature tensor } $P$  is defined by
(see, for example, [9]) components
$P_{lijk}=R_{lijk}-\frac{1}{n-1}(g_{lj}R_{ik}-g_{lk}R_{ij})$ in
arbitrary map $(U,\varphi)$ with a local coordinate system
$\{x^0,x^1,x^2,x^3\}$. For $n>2$ conversion into zero of tensor
$P$ characterises manifold of constant curvature (see [9]). Easy
to show that due to (1.3) covariant differential is
$d^*P=-\frac{n-2}{n-1}d\;Ric$. We will say that Weyl projective
curvature tensor is {\it a harmonic tensor}  if $d^*P=0$.
Definition is explained by the fact that from condition $d^*P=0$
automatically follows Bianchi identities $dP=0$. Therefore if
tensor $P$ is considered as second form
$P:\Lambda^2(TM)\rightarrow\Lambda^2(TM)$, then this form will be
simultaneously closed and coclosed and so a harmonic form (see,
for example, [13]). Condition for the Weyl projective curvature
tensor $P$  to be harmonic leads us to Codazzi equations (4.2),
which are equivalent as was proved to (4.1). Following theorem is
true:

{\bf Theorem 2.} {\it The Einstein equations belong to the class
$\Omega_2$  if and only if the Weyl projective curvature tensor
$P$ of space-time $(M,g)$ is harmonic.}

 It is well known (see, for
example, [13]) that connection $\hat{\nabla}$  in main bundle
$\pi:E\rightarrow M$ over pseudo-Riemannian manifold $(M,g)$ with
fiber metrics $g_E$ is called the {\it Young-Mills field}, if its
curvature $\hat{R}$ together with Bianchi identities $d\hat{R}=0$
satisfies the {\it Young-Mills equations} $d^*\hat{R}=0$.

 If we
consider that $E=TM$, $g_E=g$ and we take Levi-Civita connection
$\nabla$ of pseudo-Riemannian metrics $g$ as connection
$\hat{\nabla}$, then {\it Young-Mills equations} $d^*R=0$   due to
equalities (1.3) will become Codazzi equations (4.2), which are
equivalent to equations (4.1). Following theorem is true:

{\bf Theorem 3.} {\it The Einstein equations belong to the class
$\Omega_2$ if and only if the Levi-Civita connection $\nabla$ of
metrics $g$ of space-time $(M,g)$ considered as connection in
bundle $TM$ is the Young-Mills field.}

{\bf Remark.} Yang-Mills theory and the other variational theory
as Seiberg-Witten theory have been developed greatly and
influenced to topology and physics (see [14]; [15]; [16]; [17] and
etc), especially in the case of 4-dimensional manifolds.
Yang-Mills theory appeared in differential geometry (see, for
example, [18], p. 443) as pseudo-Riemannian manifolds with
harmonic curvature. This means that pseudo-Riemannian manifolds
$(M, g)$ of which curvature tensor $R$ of the Levi-Civita
connection $\nabla$ satisfies  $d^*R=0$, i.e. $\nabla$ is a
Yang-Mills connection, taking $E = TM$, the tangent bundle of $M$,
and $g_E=g$ as in our pieces 4.

\section{The class  $\Omega_3$ and geodesic mappings}

 The class $\Omega_2$  of the Einstein equations is selected
 via condition (see, [8])
$$
\nabla_kR_{ij}=\frac{1}{18}\left(4(\partial_ks)g_{ij}+(\partial_is)g_{kj}+
(\partial_js)g_{ik}\right). \eqno(5.1)
$$
Let us recall that diffeomorphism
$f:(M,g)\rightarrow(\bar{M},\bar{g})$ of pseudo-Riemannian
manifolds of dimension $n\geq 2$  is called {\it geodesic mapping}
(see [19], p. 70), if each geodesic curve of manifold $(M,g)$ is
transferred via this mapping onto geodesic curve of manifold
$(\bar{M},\bar{g})$.

 An arbitrary $(n+1)$-dimensional $(n\geq1)$  pseudo-Riemannian
manifold $(M,g)$  which scalar curvature $s\neq const$ and the
Ricci tensor $Ric$ satisfies the following equations:
$$
\nabla_kR_{ij}=\frac{n-2}{2(n-1)(n+2)}\left(\frac{2n}{n-2}(\partial_ks)g_{ij}
+(\partial_is)g_{kj}+(\partial_js)g_{ik}\right) \eqno(5.2)
$$
for a local coordinate system $\{x^0,x^1,\dots,x^n\}$  is called
{\it a manifolds $L_{n+1}$}. These manifolds were introduced by
N.S. Sinyukov (see [19], pp. 131-132). Manifold $L_{n+1}$  are an
example of pseudo-Riemannian manifolds with non-constant curvature
which admit {\it non-trivial geodesic mappings}. Easy to check
that for $n  = 3$ equations (5.1) and (5.2) are equivalent and
that is why spaces-times $(M,g)$  for which $T\in\Omega_3$ are
{\it Sinyukov spaces} $L_{n+1}$ and due to this reason admit
nontrivial geodesic mappings.

 Following theorem is true:

 {\bf Theorem 4.}  {\it The Einstein equations belong to the class
 $\Omega_2$ if and only if
the space-time $(M,g)$ is a Sinyukov space $L_4$.}

 In addition, we
call that if $(M,g)$ is an $(n+1)$-dimensional Sinyukov manifold
then the metric form $ds^2$ of the manifold $(M,g)$   has the
following form (see [19], pp. 111-117 and [20])
$$
ds^2=g_{00}(x^0)(dx^0)^2+f(x^0)\sum\limits_{a,b=1}^ng_{ab}(x^1,\dots,x^n)dx^a
dx^b\eqno(5.3)
$$
 for special local coordinate system  $\{x^0,x^1,\dots,x^n\}$. The function  $f(x^0)$,
components $g_{00}(x^0)$ and $g_{ab}(x^1,\dots,x^n)$  of the
metric form (5.3) satisfy some defined properties which can be
found in the paper [20]. From this we conclude that a space-time
$(M,g)$ with the Einstein equations of the class  $\Omega_2$ is a
warped product space-time and the Robertson-Wolker space-time, in
particular (see, for example, [21]).

\section{Three other classes}

 In brief we will describe remaining three classes of the Einstein
  equations. At the beginning lets consider class of the Einstein
  equations $\Omega_1\oplus\Omega_2$  defined by the following condition
  $g^{ij}\nabla_k T_{ij}=0$ or equivalent condition  $s=const$.
  Following theorem is true:

{\bf Theorem 5.} {\it The Einstein equations belong to the class
$\Omega_1\oplus\Omega_2$ if and only if the metric $g$  of the
space-time $(M,g)$  has constant scalar curvature.}

 Next class of
the Einstein equations $\Omega_2\oplus\Omega_3$ is defined by the
following condition  $d\left[T-\frac{1}{3}(trace_gT)g\right]=0$.
To obtain its analytical characterisation lets use the {\it Weyl
tensor} $W$ of conformal curvature (see, for example, [9]). It is
known that on arbitrary pseudo-Riemannian manifold $(M,g)$  of
dimension $n\geq3$ this tensor follows equation (see, for example,
[9] and [18])
 $$
 d^*W=-\frac{n-3}{n-2}d\left[Ric-\frac{1}{2(n-1)}s\cdot g\right].  \eqno(6.1)
 $$
We'll say that the Weyl tensor of conformal curvature is a {\it
harmonic tensor} if $d^*W=0$  (see also [18]). Definition is
explained by the fact that from condition $d^*W=0$  (see [9])
automatically follows Bianchi identities  $dW=0$. Therefore if the
Weyl tensor $W$ of conformal curvature is considered as second
form $W:\Lambda^2(TM)\rightarrow \Lambda^2(TM)$, then this form
will be simultaneously closed and coclosed and so harmonic (see
[13]).

 If the equation $T-\frac{1}{3}(trace_gT)g=Ric-\frac{1}{6}s\cdot g$ is true,
 then for space-time $(M,g)$  due to (6.1)
condition $\nabla T \in \Omega_2(M)\oplus\Omega_3(M)$ is
equivalent to harmonicity requirement $d^*W=0$ for the Weyl tensor
$W$ of conformal curvature.

Following was proved:

{\bf Theorem 6.} {\it The Einstein equations belong to the class
$\Omega_2\oplus\Omega_3$ if and only if the tensor  $W$ of
conformal curvature of the space-time $(M,g)$ is harmonic.}

 The
example of a space-time $(M,g)$  with $\nabla T\in
\Omega_2(M)\oplus \Omega_3(M)$ is a conformally flat space-time
$(M,g)$, this means that $W\equiv 0$  (see [9], p 116) and,
consequently, condition $d^*W=0$ is fulfilled automatically.

The last, third, class  $\Omega_1\oplus\Omega_3$ is defined by the
following condition \\$\delta^*\left[T-\frac{1}{6}(trace_g
T)g\right]=0$. On the basis of the identity $T-\frac{1}{6}(trace_g
T)g=Ric-\frac{1}{3}s\cdot g$, our condition turns into following
differential equations
$$
\delta^*(Ric-\frac{1}{3}s\cdot g)=0:
$$
$$
\nabla_kR_{ij}+\nabla_iR_{jk}+\nabla_jR_{ki}=\frac{1}{3}\left((\partial_ks)
g_{ij}+(\partial_is)g_{jk}+(\partial_js)g_{ki}\right). \eqno(6.2)
$$
In this case $R_{ij}\frac{d x^i}{ds}\frac{dx^j}{ds}=const$ will be
first quadratic integral of light-like geodesic equations
$x^k=x^k(s)$ of the space-time  $(M,g)$. It is obvious that on a
Sinyukov manifold $L_4$ equations (6.2) become identity.

 {\bf Remark.} The
seventh class of the Einstein equations is defined by the property
$\nabla T=0$. Chaki and Ray have studied a space-time with
covariant-constant energy-momentum tensor $T$ (see [22]). In
addition, two particular types of such space-times were considered
and the nature of each was determined (see also [23]).

 \end{document}